# Conditional-sum-of-squares estimation of models for stationary time series with long memory

P. M. Robinson[1],*

*London School of Economics*

**Abstract:** Employing recent results of Robinson (2005) we consider the asymptotic properties of conditional-sum-of-squares (CSS) estimates of parametric models for stationary time series with long memory. CSS estimation has been considered as a rival to Gaussian maximum likelihood and Whittle estimation of time series models. The latter kinds of estimate have been rigorously shown to be asymptotically normally distributed in case of long memory. However, CSS estimates, which should have the same asymptotic distributional properties under similar conditions, have not received comparable treatment: the truncation of the infinite autoregressive representation inherent in CSS estimation has been essentially ignored in proofs of asymptotic normality. Unlike in short memory models it is not straightforward to show the truncation has negligible effect.

## 1. Introduction

Consider a real-valued, strictly and covariance stationary time series $x_t$, $t \in \mathbb{Z}$. It is believed that $x_t$ has a parametric autoregressive (AR) representation

$$\sum_{j=0}^{\infty} \alpha_j(\theta_0) x_{t-j} = \varepsilon_t, \quad t \in \mathbb{Z}. \tag{1.1}$$

Here $\varepsilon_t$ is a sequence of zero-mean, uncorrelated and homoscedastic random variables, with variance $\sigma_0^2$, the $\alpha_j(\theta)$ are given functions with $p \times 1$ vector argument $\theta$, $\theta_0$ is an unknown $p \times 1$ vector, and $\alpha_0(\theta) \equiv 1$ for all $\theta$.

The range of structures $\{\alpha_j(\theta)\}$ covered by (1.1) is very broad, but of interest to us are ones which allow $x_t$ to have long memory. Usually, these are "fractional", where it is assumed that the function

$$\alpha(s;\theta) = \sum_{j=0}^{\infty} \alpha_j(\theta) s^j, \tag{1.2}$$

with complex-valued argument $s$ on the unit circle, is of form

$$\alpha(s;\theta) = (1-s)^{\delta(\theta)} \alpha^*(s;\theta), \tag{1.3}$$

where $\delta(\theta)$ is a scalar function of $\theta$ such that

$$0 < \delta(\theta_0) < \frac{1}{2} \tag{1.4}$$

[1]London School of Economics, e-mail: p.m.robinson@lse.ac.uk
*Research supported by ESRC Grant R000239936.

and

(1.5) $$0 < |\alpha^*(s;\theta_0)| < \infty, \quad |s| = 1.$$

It follows that $x_t$ has spectral density of form

(1.6) $$f(\lambda) = \frac{\sigma_0^2}{|\alpha(e^{i\lambda};\theta_0)|^2} = \sigma_0^2 \frac{|1 - e^{i\lambda}|^{-2\delta(\theta_0)}}{|\alpha^*(e^{i\lambda};\theta_0)|^2}.$$

The leading choice of $\alpha^*$ is a rational function of $s$, in which case $x_t$ is said to be a fractional autoregressive integrated moving average (FARIMA) model; $\delta(\theta_0)$ is caled the memory parameter.

Leading methods of estimation of $\theta_0$, given observations $x_1, \ldots, x_n$, are Gaussian maximum likelihood (ML), and approximations thereto. They are "approximations" in the sense that under similar conditions they have the same asymptotic normal distribution as ML, and are thus asymptotically efficient under Gaussianity. At the same time, under many departures from Gaussianity, though the efficiency is lost the limit normal distribution of all these estimates is unaffected. Assuming Gaussianity, asymptotic normality of one form of approximation, a Whittle estimate involving integration over frequency, was first established by Fox and Taqqu [4], and then by Dahlhaus [3] in case of ML estimation. Giraitis and Surgailis [5] established asymptotic normality for the estimate considered by Fox and Taqqu [4] when $\varepsilon_t$ need not be Gaussian but is independent and identically distributed with finite fourth moment. Due to the pole in the spectral density at $\lambda = 0$ (see (1.6)), the asymptotic normality proofs are considerably more challenging than those of Hannan [6] for short memory time series models, incisive though these were for such models.

An alternative estimate that has been considered in the literature is conditional-sum-of-squares (CSS) estimation, which was previously employed by Box and Jenkins [1] for short memory time series models. Define

(1.7) $$e_t(\theta) = \sum_{j=0}^{t-1} \alpha_j(\theta) x_{t-j},$$

(1.8) $$s_n(\theta) = \frac{1}{n} \sum_{t=1}^{n} e_t^2(\theta),$$

and estimate $\theta_0$ by

(1.9) $$\hat{\theta}_n = \arg\min_{\theta \in \Theta} s_n(\theta),$$

where $\Theta \subset \mathbb{R}^p$ is a compact set.

One can motivate $\hat{\theta}_n$ by the hope that $s_n(\theta_0)$ is a good approximation to $n^{-1} \times \sum_{t=1}^{n} \varepsilon_t^2$, which is itself proportional to the exponent in the density function of independent identically distributed zero-mean normal variates. Thus one hopes that (after centering at $\theta_0$ and $n^{\frac{1}{2}}$ norming) $\hat{\theta}_n$ has the same limit distributional properties as the Gaussian ML and Whittle estimates mentioned previously.

Given an initial consistency proof of $\hat{\theta}_n$, a standard approach to proving asymptotic normality entails applying the mean value theorem to $r_n(\hat{\theta}_n)$ about $\theta_0$, where

(1.10) $$r_n(\theta) = \frac{\partial s_n(\theta)}{\partial \theta} = \frac{2}{n} \sum_{t=1}^{n} h_t(\theta) e_t(\theta),$$

for

$$(1.11) \qquad h_t(\theta) = \frac{\partial e_t(\theta)}{\partial \theta}.$$

The main part of the proof then involves establishing that $n^{\frac{1}{2}} r_n(\theta_0)$ converges in distribution to a zero-mean normal vector. If the $\varepsilon_t$ are assumed to be conditionally homoscedastic martingale differences, and conditions ensuring that $h_t(\theta)$ has finite variance are imposed, such convergence is easily seen to hold (see e.g. [2]) for

$$(1.12) \qquad r_n^*(\theta_0) = \frac{2}{n} \sum_{t=1}^n h_t \varepsilon_t,$$

where $h_t = h_t(\theta_0)$. However this is only useful if also

$$(1.13) \qquad r_n^*(\theta_0) - r_n(\theta_0) = o_p\left(n^{-\frac{1}{2}}\right),$$

in other words, if the effect of replacing $e_t = e_t(\theta_0)$ by $\varepsilon_t$ is sufficiently small. Unlike the $h_t \varepsilon_t$, the $h_t e_t$ and $h_t(e_t - \varepsilon_t)$ are not zero-mean, orthogonal random variables. We can employ the Schwarz inequality:

$$(1.14) \qquad E\left|r_n^*(\theta_0) - r_n(\theta_0)\right| \leq 2n^{-1} \sum_{t=1}^n \left[E(e_t - \varepsilon_t)^2 E \|h_t(\theta_t(\theta_0))\|^2\right]^{\frac{1}{2}}.$$

Then if, say, it were true that $E(e_t - \varepsilon_t)^2 = O_p(t^{-1-\eta})$ for some $\eta > 0$, the right hand side of (1.14) would be $O_p(n^{-\frac{1}{2} - \frac{\eta}{2}})$, and (1.13) established. For short memory models $E(e_t - \varepsilon_t)^2$ typically decays fast enough, indeed even exponentially. But under quite general conditions permitting long memory (see [8]),

$$(1.15) \qquad E(e_t - \varepsilon_t)^2 \leq K t^{-1}$$

only, where $K$ is an arbitrarily large generic constant, which is insufficient to establish (1.13) using (1.14).

A more delicate proof of (1.13) is required, and this was given by Robinson [8]. As discussed there, this delicacy relates to that seen in the proofs of Fox and Taqqu [4] and others for alternative estimates of $\theta_0$. Indeed, given that these estimates and CSS should have the same limit distributional properties, it would be extraordinary if the proof for CSS were very much easier than for the other estimates.

A central limit theorem for $\hat{\theta}_n$ is given in Section 3. Prior to that, in the following section, we provide the almost sure convergence of $\hat{\theta}_n$ (under somewhat more general conditions). Hannan [6] proved this for various estimates, assuming strict stationarity and ergodicity, which is consistent with long memory. However, he did not cover CSS estimation.

## 2. Almost sure convergence

In the present section we do not require that $x_t$ necessarily has spectral density of form (1.6), with (1.5) holding, but simply that it is a zero-mean strictly stationary, ergodic process with AR representation (1.1), with the sentence after (1.1) holding, and also $\theta_0 \in \Theta$, for all $\theta \in \Theta \setminus \{\theta_0\}$

$$(2.1) \qquad \alpha(s; \theta) \neq \alpha(s; \theta_0)$$

on a subset of $|s| = 1$ of positive measure, $|\alpha(s;\theta)|$ is continuous in $\theta$ for all $s : |s| = 1$, and

$$\sum_{j=0}^{\infty} \sup_{\theta \in \Theta} |\alpha_j(\theta)| < \infty. \tag{2.2}$$

Condition (2.1) is a standard identifiability condition, and (2.2) is reasonable in that long memory models (e.g. (1.6), such as FARIMAs) typically have AR representations with summable coefficients. Note that this setup allows the spectral density to have poles at non-zero frequencies (as in certain cyclic and seasonal models), whereas (1.6) does not, in view of (1.5).

**Theorem 1.** *Under the above conditions*

$$\lim_{n \to \infty} \hat{\theta}_n = \theta_0, \quad a.s. \tag{2.3}$$

*Proof.* Theorem 1 of Hannan [6] and Theorem 1 of Fox and Taqqu [4] cover the estimate

$$\tilde{\theta}_n = \arg\min_{\Theta} s_n^{\dagger}(\theta), \tag{2.4}$$

where $s_n^{\dagger}(\theta)$ is the objective function for the integral form of Whittle estimate, i.e. $\overline{\sigma}_N^2(\theta)$ of Hannan [6] or $\sigma_N^2(\theta)$ of Fox and Taqqu [4]. We can write

$$s_n^{\dagger}(\theta) = c_n(0)\xi_0(\theta) + 2\sum_{j=1}^{n-1} c_n(j)\xi_j(\theta), \tag{2.5}$$

where

$$c_n(j) = \frac{1}{n}\sum_{t=1}^{n-j} x_t x_{t+j}, \quad 0 \le j \le n-1, \tag{2.6}$$

$$\xi_j(\theta) = \sum_{k=0}^{\infty} \alpha_k(\theta)\alpha_{k+j}(\theta). \tag{2.7}$$

From Theorem 1 of Hannan [6], and its proof, it is clear that it suffices to show that

$$\lim_{n \to \infty} \sup_{\Theta} \left| s_n^{\dagger}(\theta) - s_n(\theta) \right| = 0, \quad a.s. \tag{2.8}$$

Now

$$\begin{aligned} s_n^{\dagger}(\theta) - s_n(\theta) &= \frac{1}{n}\sum_{t=1}^{n} x_t^2 \sum_{k=n-t+1}^{\infty} \alpha_k^2(\theta) \\ &\quad + \frac{2}{n}\sum_{j=1}^{n-1}\sum_{t=1}^{n-j} x_t x_{t+j} \sum_{k=n-t-j+1}^{\infty} \alpha_k(\theta)\alpha_{k+j}(\theta) \\ &= \sum_{i=1}^{4} a_{in}(\theta), \end{aligned} \tag{2.9}$$

where

$$(2.10) \quad a_{1n}(\theta) = \gamma(0)\left\{\frac{1}{n}\sum_{j=1}^{n-1} j\alpha_j^2(\theta) + \sum_{j=n}^{\infty}\alpha_j^2(\theta)\right\},$$

$$(2.11) \quad a_{2n}(\theta) = \frac{1}{n}\sum_{t=1}^{n}\{x_t^2 - \gamma(0)\}\sum_{k=n-t+1}^{\infty}\alpha_k^2(\theta),$$

$$(2.12) \quad a_{3n}(\theta) = \frac{2}{n}\sum_{j=1}^{n-1}\gamma(j)\sum_{t=1}^{n-j}\sum_{k=n-t-j+1}^{\infty}\alpha_k(\theta)\alpha_{k+j}(\theta),$$

$$(2.13) \quad a_{4n}(\theta) = 2\sum_{j=1}^{n-1}\left\{\frac{1}{n}\sum_{t=1}^{n-j}(x_t x_{t+j} - \gamma(j))\sum_{k=n-t-j+1}^{\infty}\alpha_k(\theta)\alpha_{k+j}(\theta)\right\},$$

where

$$(2.14) \quad \gamma(j) = \operatorname{cov}(x_0, x_j).$$

It remains to prove

$$(2.15) \quad \lim_{n\to\infty}\sup_{\theta\in\Theta}|a_{in}(\theta)| = 0 \quad \text{a.s.,} \quad i = 1, 2, 3, 4.$$

As the proofs for $i = 1, 2$ are similar to but simpler than those for $i = 3, 4$, we give only the latter. We have

$$(2.16) \quad \sup_{\Theta}|a_{3n}(\theta)| \leq \frac{2}{n}\sum_{j=1}^{n-1}|\gamma(j)|\left\{\sum_{j=0}^{\infty}\sup_{\theta\in\Theta}|\alpha_j(\theta)|\right\}^2.$$

The quantity in braces is finite and since, by the Riemann-Lebesgue theorem, existence of a spectral density implies $\lim_{j\to\infty}\gamma(j) = 0$, it follows from the Toeplitz lemma that $(2.16) \to 0$ as $n \to \infty$. Next, by summation-by-parts

$$(2.17) \quad \begin{aligned} a_{4n}(\theta) = &-2\sum_{j=1}^{n-1}\sum_{t=1}^{n-j-1}\frac{t}{n}\{c_t(j) - \gamma(j)\}\alpha_{n-t-j+1}(\theta)\alpha_{n-t+1}(\theta) \\ &+ 2\sum_{j=1}^{n-1}\frac{1}{n}\sum_{t=1}^{n-j}\{x_t x_{t+j} - \gamma(j)\}\sum_{k=1}^{\infty}\alpha_k(\theta)\alpha_{k+j}(\theta). \end{aligned}$$

The modulus of the first term on the right has supremum, over $\Theta$, bounded by

$$(2.18) \quad K\sum_{t=1}^{n}\sup_{1\leq j\leq n}|c_t(j) - \gamma(j)|\sup_{\Theta}|\alpha_{n-t+1}(\theta)|$$

using (2.2). Using (2.2) again, and Theorem 1 of Hannan [7] and the Toeplitz lemma, it follows that (2.18) is $o(1)$ a.s. The second term in (2.17) can be similarly handled. □

## 3. Asymptotic normality

We assume now in addition that $x_t$ has spectral density (1.6), with (1.4), (1.5) satisfied, that $\theta_0$ is an interior point of $\Theta$, that the $\varepsilon_t$ in (1.1) are independent with zero mean, variance $\sigma_0^2$ and uniformly bounded fourth moment, that $\alpha(s;\theta)$ is twice continuously differentiable in $\theta$, and that the matrix

$$(3.1) \qquad \Omega = \frac{1}{2\pi}\int_{-\pi}^{\pi}\left[\begin{array}{c}\log\left|1-e^{i\lambda}\right|^2 \\ -2\frac{\partial}{\partial\theta}\log\left|\alpha\left(e^{i\lambda};\theta_0\right)\right|\end{array}\right]\left[\begin{array}{c}\log\left|1-e^{i\lambda}\right|^2 \\ -2\frac{\partial}{\partial\theta}\log\left|\alpha\left(e^{i\lambda};\theta_0\right)\right|\end{array}\right]' d\lambda$$

is positive definite.

**Theorem 2.** *Under the above conditions, as $n\to\infty$ $n^{\frac{1}{2}}(\hat{\theta}_n - \theta_0)$ converges in distribution to a p-variate normal vector with zero mean and covariance matrix $\Omega^{-1}$.*

*Proof.* As discussed in Section 1, we have

$$(3.2) \qquad 0 = r_n(\hat{\theta}_n) = r_n(\theta_0) + \tilde{T}_n(\hat{\theta}_n - \theta_0),$$

where $\tilde{T}_n$ is the matrix formed by evaluating, for $i=1,\ldots,p$, the $i$-th row of the matrix $T_n(\theta) = (\partial^2/\partial\theta\partial\theta')s_n(\theta)$ at $\theta = \tilde{\theta}_i$, where $\|\tilde{\theta}_i - \theta_0\| \leq \|\hat{\theta}_n - \theta_0\|$, $\|\cdot\|$ denoting Euclidean norm.

Define

$$(3.3) \qquad \zeta_j = \frac{\partial}{\partial\theta}\alpha_j(\lambda;\theta),$$

so that

$$(3.4) \qquad h_t = \sum_{j=1}^{t-1}\zeta_j x_{t-j},$$

and define also

$$(3.5) \qquad \rho_t = \sum_{j=1}^{\infty}\zeta_j x_{t-j},$$

$$(3.6) \qquad r_n = \frac{1}{n}\sum_{t=1}^{n}\rho_t\varepsilon_t.$$

Write $r_n(\theta_0) - r_n = r_{1n} + r_{2n} + r_{3n}$, where

$$(3.7) \qquad r_{1n} = 2n^{-1}\sum_{t=1}^{n}(h_t - \rho_t)\varepsilon_t,$$

$$(3.8) \qquad r_{2n} = 2n^{-1}\sum_{t=1}^{n}\rho_t(e_t - \varepsilon_t),$$

$$(3.9) \qquad r_{3n} = 2n^{-1}\sum_{t=1}^{n}(h_t - \rho_t)(e_t - \varepsilon_t).$$

We show that $r_{in} = o_p(n^{-\frac{1}{2}})$, $i = 1,2,3$. To deal with $r_{1n}$, we may write

$$(3.10) \qquad h_t - \rho_t = -\sum_{j=t}^{\infty}\zeta_j x_{t-j} = -\sum_{j=1}^{\infty}\chi_{jt}\varepsilon_{-j},$$

where

$$\chi_{jt} = \sum_{k=0}^{j} \zeta_{k+j}\beta_{j-k}. \tag{3.11}$$

Since

$$E\|h_t - \rho_t\|^2 = \sigma_0^2 \sum_{j=1}^{\infty} \|\chi_{jt}\|^2 \leq K\frac{(\log t)^2}{t} \tag{3.12}$$

as noted on p. 1824 of [8], and $\varepsilon_t$ is independent of $h_t - \rho_t$, it follows that

$$E\|r_{1n}\|^2 \leq \frac{K}{n^2}\sum_{t=1}^{n} t^{-1} \leq K\frac{\log n}{n^2}. \tag{3.13}$$

Next, we can write

$$e_t - \varepsilon_t = -\sum_{j=1}^{\infty} \lambda_{jt}\varepsilon_{-j}, \tag{3.14}$$

where

$$\lambda_{jt} = \sum_{k=0}^{j} \alpha_{k+j}\beta_{t-k}. \tag{3.15}$$

Thus, from Lemma 16 of Robinson [8],

$$E\|r_{2n}\|^2 \leq K\frac{(\log n)^3}{n^2}. \tag{3.16}$$

Finally,

$$E\|r_{3n}\| \leq \frac{1}{n}\sum_{t=1}^{n}\left(E\|h_t - \rho_t\|^2 E(e_t - \varepsilon_t)^2\right)^{\frac{1}{2}}$$
$$\leq \frac{K}{n}\sum_{t=1}^{n}\frac{\log t}{t}$$
$$\leq K\frac{(\log n)^2}{n}, \tag{3.17}$$

using (3.12) and also Lemma 14 of Robinson [8]. This completes the proof that $r_{in} = o_p(n^{-\frac{1}{2}})$, $i = 1, 2, 3$. The remainder of the proof is easier, and more standard, and is omitted. □

**Acknowledgment**

I thank a referee for a careful reading of the paper.